\documentclass[12pt,reqno]{amsart}
\usepackage{mathrsfs}
\usepackage{amsfonts}
\usepackage{amsfonts}
\usepackage{amsfonts}
\usepackage{amsfonts}
\usepackage{amsfonts}
\usepackage{amsfonts}
\usepackage{amsfonts}
\usepackage{amsfonts}
\usepackage{amsfonts}
\usepackage{amsfonts}
\usepackage{amsfonts}
\usepackage{amsfonts}
\usepackage{amsfonts}
\usepackage{amsfonts}
\usepackage{amsfonts}
\usepackage{amsfonts}
\usepackage{amsfonts}
\usepackage{amsfonts}
\usepackage{a4wide}
\let\Section=\section
\def\section{\setcounter{equation}{0}\Section}

  \sf
\begin{document}
\title{}
\begin{center}
 {\Large\bf  Carleson measures  for Besov-Sobolev  spaces }
\end{center}
\begin{center}
 {\Large\bf   with applications in the unit ball of ${\bf
C}^{n*}$ }
\end{center}
\vspace{0.5cm}
\author{}
 \begin{center}
{\small
 Ru  Peng$^{1}$\ \ \ \ \ \ \ \ Caiheng  Ouyang$^{2}$}
\end{center}
\vspace{0.2cm}
\begin{center}

{\footnotesize 1.\  Department of Mathematics, Wuhan University of Technology, Wuhan 430070, China}\\
{\footnotesize 2.\  Wuhan Institute of Physics and Mathematics, The Chinese Academy of Sciences, Wuhan 430071, China}\\
\end{center}
\vspace{0.4cm}
\thanks{*Supported in part by the National
Natural Science Foundation of China (Nos. 10971219, 11126048 and
11101279) and the Fundamental Research Funds for the Central
Universities (No. 2012-Ia-018).}

{\footnotesize\bf Abstract.} {\footnotesize This paper is devoted to
give the connections between Carleson measures for Besov-Sobolev
spaces $B_p^\sigma (B)$  and $p$-Carleson measure in the unit ball
of ${\bf C}^n$. As applications, we  characterize the
Riemann-Stieltjes operators and multipliers acting on $B_p^\sigma
(B)$ spaces by means of Carleson measures for $B_p^\sigma (B)$. }
\vspace{0.1cm}

\vspace{0.1cm}


\date{}\maketitle{}

\vskip 0.1cm \arraycolsep1.5pt
\newtheorem{Lemma}{Lemma}[section]
\newtheorem{Theorem}{Theorem}[section]
\newtheorem{Definition}{Definition}[section]
\newtheorem{Proposition}{Proposition}[section]
\newtheorem{Remark}{Remark}[section]
\newtheorem{Corollary}{Corollary}[section]

\small

\begin{center}
{\small\bf\S1\ Introduction}
\end{center}
\vspace{0.3cm}

Let $B=\{{z\in {\bf C}^n:\mid z\mid< 1}\}$ be the unit ball of ${\bf
C}^n$ $(n>1)$, $S=\{z\in {\bf C}^n: \mid z\mid=1\}$ be its boundary.
Let $d\upsilon$ denote the normalized Lebesgue measure of $B$, i.e.
$\upsilon (B)=1$, and $d\sigma$ denote the normalized rotation
invariant Lebesgue measure of $S$ satisfying $\sigma (S)=1$. Let
$d\lambda (z)={(1-{\mid z\mid}^2)}^{-n-1}d\upsilon (z)$ be the
invariant measure on the ball.

We denote the class of all holomorphic functions in $B$ by $H(B)$.
In \cite{ARS2}, for integer $m> 0$, and for $0\leq \sigma <\infty$,
$1< p< \infty$, $m+\sigma >{n\over p}$,  the  Besov-Sobolev spaces
$B_p^\sigma (B)$
 are defined to consist of those  $f\in H(B)$ on the
ball such that
$${\left\{\sum\limits_{k=0}^{m-1}{|f^{(k)}(0)|}^p+\int_B {|{(1-{|z|}^2)}^{m+\sigma}f^{(m)}(z)|}^p d \lambda (z)\right\}}^{1\over p}<\infty.
\ \ \ \ \ \ \ \ \ \ \ \ \ \ \ \ (1.1)$$ Here $f^{(m)}$ is the
$m^{th}$ order complex derivative of $f$. The spaces $B_p^\sigma
(B)$ are independent of $m$ and are Banach spaces with norms given
in (1.1).

For $p=2$, these are Hilbert spaces with the usual inner product in
${\bf C}^n$. This scale of spaces includes the Dirichlet spaces
$B_2(B)=B_2^0(B)$, weighted Dirichlet-type spaces with
$0<\sigma<{1\over 2}$, the Drury-Arveson Hardy spaces
$H_n^2=B_2^{1\over 2}(B)$, the Hardy spaces $H^2=B_2^{n\over 2}
(B)$, and the weighted Bergman spaces with $\sigma>{n\over 2}$ (see
\cite{OYZ1} and \cite{Zh2}).

For $f\in H(B), z\in B$, its complex gradient and invariant gradient
are defined as
\begin{center}
$\nabla f(z)={\nabla}_zf=({\frac{\partial f}{\partial
z_1}(z)},\ldots,{\frac{\partial f}{\partial z_n}(z)})$,\ \ \ \
$\widetilde{\nabla}f(z)=\nabla(f\circ\varphi_z)(0)$,
\end{center}
where $\varphi _z$ is the M${\rm \ddot{o}}$bius transformation for
$z\in B$ , which satisfies $\varphi _z (0)=z$, $\varphi _z(z)=0$ and
$\varphi _z\circ\varphi _z=I$,
 and its radial derivative $Rf(z)=<\nabla f(z),\bar{z}>=\sum\limits_{j=1}^{n}\frac{\partial f}{\partial
 z_j}(z)z_j$. In \cite{Zh2}, the invertible "radial" operators
 $R^{\alpha,t}: H(B)\rightarrow H(B)$ is denoted by
 $$R^{\alpha,t}f(z)=\sum\limits_{k=0}^{\infty}\frac{\Gamma (n+1+\alpha)\Gamma (n+1+k+\alpha+t)}{\Gamma (n+1+\alpha+t)\Gamma
 (n+1+k+\alpha)}f_k (z),$$
 provided neither $n+\alpha$ nor  $n+\alpha+t$ is a negative
 integer, and where $f(z)=\sum_{k=0}^{\infty} f_k (z)$ is the
 homogeneous expansion of $f$. If the inverse of $R^{\alpha,t}$ is
 denoted by $R_{\alpha,t}$, then Proposition 1.14 of \cite{Zh2}
 yields
 $$R^{\alpha,t}\left(\frac{1}{{(1-\langle
 z,w\rangle)}^{n+1+\alpha}}\right)=\frac{1}{{(1-\langle
 z,w\rangle)}^{n+1+\alpha+t}},$$
$$R_{\alpha,t}\left(\frac{1}{{(1-\langle
 z,w\rangle)}^{n+1+\alpha+t}}\right)=\frac{1}{{(1-\langle
 z,w\rangle)}^{n+1+\alpha}},$$
 for all $w\in B$. Thus for any $\alpha$, $R^{\alpha,t}$ is
 approximately differentiation of order $t$.

Using the similar method of Lemma 6.3, Theorem 6.1 and Theorem 6.4
of \cite{Zh2}, we know the definition  (1.1) is equivalent to the
seminorm
$${\| f\|}_{B_p^\sigma (B)}^p=\int_B {| R^{\alpha, m} f(z)|}^p{(1-{|z|}^2)}^{mp+\sigma p-n-1} d
v(z)<\infty,\ \ \ \ \ \ \ \ \ \ (1.2)$$ for integer $m> 0$, and for
$0\leq \sigma <\infty$, $1< p< \infty$, $m+\sigma >{n\over p}$,
where neither $n+\alpha$ nor $n+\alpha+m$ is a negative integer.

A positive Borel measure $\mu$ on $B$ is called a Carleson measure
for $B_p^\sigma (B)$  if there is a constant $C>0$ such that
$$\int_B {| f(z)|}^p d\mu (z)\leq C{\| f\|}_{B_p^\sigma (B)}^p,\ \ \ \
\ \ \ \forall f\in B_p^\sigma (B).$$

For $z\in B$ and $r>0$, we  denote  $E(z,r)=\{w\in
B:|\varphi_z(w)|<r\}$ the pseudo-hyperbolic metric ball at $z$. For
$\xi\in S$ and $\delta >0$, let $Q_\delta(\xi)=\{z\in B:|1-\langle
z,\xi\rangle|<\delta\}$. For a positive Borel measure $\mu$ on $B$,
if
$${\|\mu\|}_{CM_p}^2=\sup\{{\mu(Q_\delta(\xi))\over \delta^{np}};\xi\in S, \delta>0\}<\infty,$$ we call $\mu$ a p-Carleson measure.

The study of Carleson measures for Besov-Sobolev spaces has a long
history. It plays the important role in function spaces and operator
theory. In one variable, various authors give their
characterizations by using appropriate capacities (see \cite{KS},
\cite{St} and \cite{Wu}). Recently, N. Arcozzi, R. Rochberg and E.
Sawyer extend themselves earlier characterization in \cite{ARS} to
higher dimensions. In \cite{ARS1}, they described the Carleson
measures for $B_p^\sigma (B)$ on the unit ball in ${\bf C}^n$  for
$\sigma=0$ and $1<p<2+{1\over n-1}$ (the difficult range $p\in
[2+{1\over n-1}, \infty )$ remains open) in terms of a discrete tree
condition on the associated Bergman tree. Subsequently, in
\cite{ARS2}, they considered $0\leq\sigma\leq {1\over 2}$, and
focused their attention on the Hilbert spaces $p=2$ (the range
${1\over 2}<\sigma< {n\over 2}$ remains mysterious ,  recently, the
difficult range ${1\over 2}<\sigma< {n\over 2}$ was given in
\cite{Vo} by a ``T1 Condition", but only for $p=2$). Our starting
point is an attempt to get easier conditions to characterize
Carleson measures for $B_p^\sigma (B)$ in the unit ball of ${\bf
C}^n$ for the whole values of $p$ and $\sigma$.

In this paper, we give the connections between Carleson measures for
$B_p^\sigma (B)$ and $p$-Carleson measure in Theorem 2.1, which
seems easier to be verified than ``T1 Condition" and the discrete
tree conditions  . In Theorem 2.1, we consider not only the case
$\sigma=0$ but also the case $0<\sigma<\infty$, and these results
hold for all the ranges $1<p<\infty$. The reason  why there is a
difference of any small $\varepsilon>0$ between necessary condition
and sufficient condition is that $p$-Carleson measure is weaker than
Carleson measures for $B_p^\sigma (B)$, which was mentioned in
\cite{ARS} for the case of the unit disc. Therefore, even in
one-dimensional situation, in \cite{GP}, unified necessary and
sufficient condition holds only for $0<p<q<\infty$. This paper may
be regarded as an extension of \cite{GP} to the case $p=q$ in the
higher dimensions. But the traditional method of one complex
variable in \cite{GP} is not applicable to several complex
variables. The key point in the proof of Theorem 2.1 is that duality
theorem is adapted. At the same time, we apply these results to
characterize Riemann-Stieltjes operators and multipliers for
$B_p^\sigma (B)$ in the unit ball of ${\bf C}^n$.

 The Riemann-Stieltjes operators $V_\varphi$ and $U_\varphi$  with the holomorphic
symbol $\varphi$ on B are defined as follows (see \cite{Hu},
\cite{Xi1}) :
$$V_\varphi f(z)=\int_0^1 f(tz)R\varphi (tz)\frac{dt}{t},\ \ \ \ \ \
\ U_\varphi f(z)=\int_0^1 \varphi(tz)Rf (tz)\frac{dt}{t},\ \ \ z\in
B.$$

It is easy to see that the pointwise multipliers  $M_\varphi$  are
determined by

$$M_\varphi f(z)=\varphi (z)f(z)=\varphi (0)f(0)+V_\varphi f(z)+U_\varphi f(z),\ \ \ \ \ z\in
B.$$ Of course, in the above definition $f$ is assumed to be
holomorphic in $B$. Clearly, $V_\varphi f=U_f \varphi$ and the
Riemann-Stieltjes operator can be viewed as a generalization of the
well known Ces$\acute{a}$ro operator.

Throughout this paper, $C$, $M$ denote positive constants which are
not necessarily the same at each appearance. The expression
$A\approx B$ means that there exists a positive $C$ such that
$C^{-1}B\leq A\leq CB$.

\vspace{0.5cm}

\begin{center}
{\small\bf\S2\ \  Carleson measures for Besov- Sobolev spaces }
\end{center}
\vspace{0.3cm} \setcounter{section}{2} \setcounter{equation}{0}

Similar to the proof of Lemma 3.2 of \cite{OF4}, it is easy to prove
the following Lemma 2.1 and Lemma 2.2. For the convenience of
readers, we give the details of the proof of Lemma 2.1.

\vspace{0.1cm}

{\bf Lemma 2.1}\ \ Let $1<p<\infty$, $\mu$ be a positive Borel
measure. Then the following statements are equivalent :

(i)\  The measure $\mu$ satisfies
$$\sup\{\mu({Q_\delta(\xi)});\xi\in
S\}\leq C \delta^{np}.$$

(ii)\ For every $s>0$,
$$\sup\{\int_B\frac{{(1-{|w|}^2)}^s}{{|1-\langle
z,w\rangle|}^{np+s}}d\mu(z);w\in B\}<\infty.$$

(iii)\ For some $s>0$,
$$\sup\{\int_B\frac{{(1-{|w|}^2)}^s}{{|1-\langle
z,w\rangle|}^{np+s}}d\mu(z);w\in B\}<\infty.$$

{\bf Proof} \ To show that (i) implies (ii). Clearly, it is
sufficient to prove (ii) for $w\in B$ near to the boundary. Let
$J_w$ be the integer part of $\left(\log_2\frac
{1}{1-{|w|}^2}\right)-1$. For $j=0,1,\dots,J_w$, consider the sets

\begin{center}

$\Omega_0=\emptyset$,\ \ \ \  $\Omega_j=\{z\in B: |1-\langle
z,w\rangle |\leq 2^{j}(1-{|w|}^2)\},j\geq 1$

\end{center}

Thus, (ii) follows from
\begin{eqnarray*}
&&\int_B\frac{{(1-{|w|}^2)}^s}{{|1-\langle
z,w\rangle|}^{np+s}}d\mu(z)\\
&\leq&C\sum\limits_{j=1}^{J_w}\frac{{(1-{|w|}^2)}^s}{2^{(np+s)(j-1)}{(1-{|w|}^2)}^{np+s}}\mu{(\Omega_j\backslash\Omega_{j-1})}\\
&&+\int_{B\backslash\Omega_{J_w}}\frac{{(1-{|w|}^2)}^s}{{|1-\langle
z,w\rangle|}^{np+s}}d\mu(z)\\
&\leq&C\sum\limits_{j=1}^{J_w}\frac{1}{2^{(np+s)(j-1)}{(1-{|w|}^2)}^{np}}\mu{(\Omega_j)}+C\mu(B\backslash\Omega_{J_w})\\
&\leq&C\sum\limits_{j=1}^{J_w}\frac{1}{2^{(np+s)(j-1)}{(1-{|w|}^2)}^{np}}{{\left(2^j
(1-{|w|}^2)\right)}^{np}}+C\\
&\leq&C\sum\limits_{j=0}^{\infty}\frac{1}{2^{sj}}<\infty.
\end{eqnarray*}

That (ii) implies (iii) is trivial. To obtain (i) from (iii), note
that for $0<\delta<1$ and $z\in {Q_\delta(\xi)}$ we have $|1-\langle
z,(1-\delta)\xi\rangle|\approx \delta$. Therefore,

\begin{eqnarray*}
\mu({Q_\delta(\xi)})&\leq
&C\delta^{np}\int_{Q_\delta(\xi)}\frac{{(1-{|(1-\delta)\xi|}^2)}^s}{{|1-\langle
z,(1-\delta)\xi\rangle|}^{np+s}}d\mu(z)\\
&\leq& C\delta^{np}.
\end{eqnarray*}

\vspace{0.1cm}

 {\bf Lemma 2.2}\ \ Let $1<p<\infty$, $\varepsilon>0$ and $\mu$ be a
positive Borel measure. Then the following statements are equivalent
:

(i)\  The measure $\mu$ satisfies
$$\sup\{\mu({Q_\delta(\xi)});\xi\in
S\}\leq C{\log^{1-p-\varepsilon}{2\over \delta}}.$$

(ii)\ For every $s>0$,
$$ \sup\{\log^{p-1+\varepsilon}{2\over
{1-{|w|}^2}}\int_B\frac{{(1-{|w|}^2)}^s}{{|1-\langle
z,w\rangle|}^{s}}d\mu(z);w\in B\}<\infty.$$

(iii)\ For some $s>0$,
$$\sup\{\log^{p-1+\varepsilon}{2\over
{1-{|w|}^2}}\int_B\frac{{(1-{|w|}^2)}^s}{{|1-\langle
z,w\rangle|}^{s}}d\mu(z);w\in B\}<\infty.$$

\vspace{0.1cm}

 {\bf Lemma 2.3}\ \ Suppose integer $m>0$,  $1<p<\infty$, $\varepsilon>0$, neither $n+\alpha$ nor $n+\alpha+m$ is a
negative integer, $\mu$ be a positive Borel measure such that
$$\sup\{{\mu(Q_\delta(\xi))\over {{\log^{1-p-\varepsilon}{2\over \delta} }}};\xi\in S, \delta>0\}<\infty.$$
Then, for $M> m p$,
$$\int_{B}{\left(\int_B | R^{\alpha,m} f(w)|\frac{{(1-{|w|}^2)}^M}{{|1-\langle
z,w\rangle|}^{n+M+1-m}}dv(w)\right)}^p d\mu(z)\leq C{\| f\|}_{B_p^0
(B)}^p.$$

{\bf Proof} \ \ Let
$$I={\left(\int_{B}{\left(\int_B |R^{\alpha,m} f(w)|\frac{{(1-{|w|}^2)}^M}{{|1-\langle
z,w\rangle|}^{n+M+1-m}}dv(w)\right)}^p d\mu(z)\right)}^{1\over p}.$$

Let $\|\cdot\|_{L_p}$ denote the usual norm on $L^p(B,d\mu)$,
${1\over p}+{1\over q}=1$. By duality,
$$I=\sup\limits_{{\|\psi\|}_{L_q}=1}\left\{\int_{B}\int_B\frac{|R^{\alpha,m} f(w)|{(1-{|w|}^2)}^M}{{|1-\langle
z,w\rangle|}^{n+M+1-m}}dv(w)|\psi(z)| d\mu(z)\right\}.$$

 By H${\rm\ddot{o}}$lder's inequality and Fubini's theorem, we have
\begin{eqnarray*}
I&\leq&\sup\limits_{{\|\psi\|}_{L_q}=1}{\left(\int_{B}\int_{B}\frac{{|
R^{\alpha,m} f(w)|}^p{(1-{|w|}^2)}^M\log^{p-1+\varepsilon}{2\over
{1-{|w|}^2}}}{{|1-\langle z,w\rangle|}^{M-mp+n+1}}dv(w)
d\mu(z)\right)}^{1\over p}\\
&\times&{\left(\int_{B}\int_{B}\frac{{|\psi(z)|}^q{(1-{|w|}^2)}^M}{{|1-\langle
z,w\rangle|}^{n+1+M}\log^{1+{q\over p}\varepsilon}{2\over
{1-{|w|}^2}}}dv(w)
d\mu(z)\right)}^{1\over q}\\
&\leq&
C\sup\limits_{{\|\psi\|}_{L_q}=1}{\left(\int_{B}\int_{B}\frac{{(1-{|w|}^2)}^{M-mp}\log^{p-1+\varepsilon}{2\over
{1-{|w|}^2}}}{{|1-\langle
z,w\rangle|}^{M-mp}} d\mu(z){| R^{\alpha,m} f(w)|}^p {(1-{|w|}^2)}^{mp-n-1} dv(w)\right)}^{1\over p}\\
&\times&{\left(\int_{B}\int_{B}\frac{{(1-{|w|}^2)}^M}{{|1-\langle
z,w\rangle|}^{n+1+M}\log^{1+{q\over p}\varepsilon}{2\over
{1-{|w|}^2}}}dv(w){|\psi(z)|}^q d\mu(z)\right)}^{1\over q}.
\end{eqnarray*}

Similar to the proof of Lemma 3.4 in \cite{OF1}, it is clear that
the inner integral of the last line above is bounded. And noting
$M>mp$, by Lemma 2.2  we can get
\begin{eqnarray*}
I&\leq&C\sup\limits_{{\|\psi\|}_{L_q}=1}{\left(\int_{B}{|
R^{\alpha,m} f(w)|}^p{(1-{|w|}^2)}^{mp-n-1}dv(w)\right)}^{1\over p}\times{\left(\int_{B}{|\psi(z)|}^q d\mu(z)\right)}^{1\over q}\\
&\leq&C{\left(\int_{B}{| R^{\alpha,m}
f(w)|}^p{(1-{|w|}^2)}^{mp-n-1}dv(w)\right)}^{1\over p}\\
&\leq&C{\| f\|}_{B_p^0 (B)}.
\end{eqnarray*}

\vspace{0.1cm}

{\bf Theorem 2.1}\ \ Suppose integer $m>0$, $0\leq \sigma<\infty$,
$1<p<\infty$, $m+\sigma> {n\over p}$, neither $n+\alpha$ nor
$n+\alpha+m$ is a negative integer. Let $\mu$ be a positive Borel
measure in $B$.

(i)\ For $\sigma>0$, if $\mu$ is a Carleson measure for $B_p^\sigma
(B)$, then $\mu$ is a ${\sigma p\over n}$-Carleson measure,
conversely, if $\varepsilon >0$ and $\mu$ is a ${\sigma
p+\varepsilon\over n}$-Carleson measure, then $\mu$ is a Carleson
measure for $B_p^\sigma (B)$.

(ii)\ For $\sigma=0$, if $\mu$ is a Carleson measure for $B_p^0
(B)$, then $\mu$ satisfy
$$\sup\{{\mu(Q_\delta(\xi))\over {{\log^{1-p}{2\over
\delta} }}};\xi\in S, \delta>0\}<\infty,$$ conversely, if
$\varepsilon>0$ and $\mu$ satisfy $$\sup\{{\mu(Q_\delta(\xi))\over
{{\log^{1-p-\varepsilon}{2\over \delta} }}};\xi\in S,
\delta>0\}<\infty,$$ then $\mu$ is a Carleson measure for $B_p^0
(B)$.

{\bf Proof}

(i)\ For $\sigma>0$, suppose first that $\mu$ is a Carleson measure
for $B_p^\sigma (B)$, then
\begin{equation}\label{2.1}
{\|f\|}_{L^p(d\mu)}\leq C{\|f\|}_{B_p^\sigma (B)}.
\end{equation}

We can find a constant $s>0$ such that $\sigma+{s\over
p}-n-1=\alpha+N$ for some positive integer $N$. Applying (2.1) to
the test functions $$f_w(z)= \frac{{(1-{|w|}^2)}^{s\over
p}}{{(1-\langle z,w\rangle)}^{\sigma+{s\over p}}}, \ \ \ \ w\in B,$$
by Lemma 2.18 of \cite{Zh2}, we can get
\begin{eqnarray*}
\int_B\frac{{(1-{|w|}^2)}^{s}}{{| 1-\langle z,w\rangle|}^{\sigma p
+s}}d\mu(z)&\leq& C
\int_B{|R^{\alpha,m} f_w (z)|}^p{(1-{|z|}^2)}^{p(\sigma+m)-n-1}dv(z)\\
&\leq& C \int_B\frac{{(1-{|w|}^2)}^s}{{|
1-\langle z,w\rangle|}^{\sigma p+mp+s}}{(1-{|z|}^2)}^{p(\sigma+m)-n-1}dv(z)\\
&=&C\int_B\frac{{(1-{|w|}^2)}^{s}{(1-{|z|}^2)}^{p(\sigma+m)-n-1}}{{|
1-\langle z,w\rangle|}^{(n+1)+[p(\sigma+m)-n-1]+s}}dv(z)\\
&\leq&C{(1-{|w|}^2)}^{s}{(1-{|w|}^2)}^{-s}\\
&\leq&C.\ \ \ \ \ \ \ \ \ \ \ \ \ \ \ \ \ \ \ \ \ \ \ \ \ \ \ \ \ \
\ (2.2)
\end{eqnarray*}
where Proposition 1.4.10 of \cite{Ru} is used. Taking
$\sup\limits_{w\in B}$ of (2.2) and  by Lemma 2.1, we know that
$\mu$ is a ${\sigma p \over n}$-Carleson measure.

On the other hand, suppose $\mu$ is a ${\sigma p+\varepsilon\over
n}$-Carleson measure, we need to prove $\mu$ is a Carleson measure
for $B_p^\sigma (B)$.

Fix a sufficiently large positive integer $K$ and let $M=\alpha+K$.
Then
$$R^{\alpha,m} f(z)=C_M\int_B R^{\alpha,m} f(w){{(1-{|w|}^2)}^M\over
{(1-\langle z,w\rangle)}^{n+1+M}}dv(w).$$ Acting on the above
equation by the inverse operator $R_{\alpha,m}$,
$$f(z)=C_M R_{\alpha,m}\int_B R^{\alpha,m} f(w){{(1-{|w|}^2)}^M\over
{(1-\langle z,w\rangle)}^{n+1+M}}dv(w).$$ By Lemma 2.18 of
\cite{Zh2}, there exists a polynomial $P(z,w)$ such that
$$f(z)=C_M \int_B {P( z,w) R^{\alpha,m} f(w){(1-{|w|}^2)}^M\over
{(1-\langle z,w\rangle)}^{n+1+M-m}}dv(w),$$
and consequently, we can
get
$$|f(z)|\leq C\int_B|R^{\alpha,m}  f(w)|{{(1-{|w|}^2)}^M\over {|1-\langle
z,w\rangle|}^{n+M+1-m}}dv(w).\ \ \ \ \ \ \ \ \ \ \ \ \ \ \ \ \ \ \ \
\ \ \ \ \ \ \ \ \ \ (2.3)$$ By (2.3) and a process similar to the
proof of Lemma 2.3, we have
\begin{eqnarray*}
& &{\left(\int_B{|f(z)|}^p d\mu (z)\right)}^{1\over p}\\
&\leq&C{\left(\int_B{\left(\int_B|R^{\alpha,m}
f(w)|{{(1-{|w|}^2)}^M\over {|1-\langle
z,w\rangle|}^{n+M+1-m}}dv(w)\right)}^p d\mu (z)\right)}^{1\over p}\\
&=&C\sup\limits_{{\|\psi\|}_{L_q}=1}\left\{\int_{B}\int_B\frac{|R^{\alpha,m}
f(w)|{(1-{|w|}^2)}^M}{{|1-\langle
z,w\rangle|}^{n+M+1-m}}dv(w)|\psi(z)| d\mu(z)\right\}\\
&\leq&C{\left(\int_{B}\int_{B}\frac{{| R^{\alpha,m}
f(w)|}^p{(1-{|w|}^2)}^M}{{|1-\langle
z,w\rangle|}^{M-mp+n+1+\varepsilon}}dv(w)
d\mu(z)\right)}^{1\over p}\\
&\times&\sup\limits_{{\|\psi\|}_{L_q}=1}{\left(\int_{B}\int_{B}\frac{{|\psi(z)|}^q{(1-{|w|}^2)}^M}{{|1-\langle
z,w\rangle|}^{n+1+M-{q\over p}\varepsilon}}dv(w)
d\mu(z)\right)}^{1\over q}\\
&\leq&
C{\left(\int_{B}\left(\int_{B}\frac{{(1-{|w|}^2)}^{M-[(\sigma+m)p-n-1]}}{{|1-\langle
z,w\rangle|}^{M-[(\sigma+m)p-n-1]+(\sigma p+\varepsilon)}} d\mu(z)\right){| R^{\alpha,m} f(w)|}^p {(1-{|w|}^2)}^{(\sigma+m)p-n-1} dv(w)\right)}^{1\over p}\\
&\times&\sup\limits_{{\|\psi\|}_{L_q}=1}{\left(\int_{B}\left(\int_{B}\frac{{(1-{|w|}^2)}^M}{{|1-\langle
z,w\rangle|}^{n+1+M-{q\over p}\varepsilon}}dv(w)\right){|\psi(z)|}^q
d\mu(z)\right)}^{1\over q}.\ \ \ \ \ \ \ \ \ \ \ \ (2.4)
\end{eqnarray*}
Applying Lemma 2.1 and Proposition 1.4.10 of \cite{Ru} to the two
inner integrals in the end of (2.4) respectively, we know
\begin{eqnarray*}
(2.4)&\leq&C{\left(\int_{B}{|
R^{\alpha,m} f(w)|}^p{(1-{|w|}^2)}^{(\sigma+m)p-n-1}dv(w)\right)}^{1\over p}\times\sup\limits_{{\|\psi\|}_{L_q}=1}{\left(\int_{B}{|\psi(z)|}^q d\mu(z)\right)}^{1\over q}\\
\\
&\leq&C {\left(\int_{B}{|
R^{\alpha,m} f(w)|}^p{(1-{|w|}^2)}^{(\sigma+m)p-n-1}dv(w)\right)}^{1\over p}\\
 &\leq&C{\| f\|}_{B_p^\sigma (B)}.
\end{eqnarray*}
this implies $\mu$ is a Carleson measure for $B_p^\sigma (B)$.

(ii)\ For $\sigma=0$, suppose $\mu$ is a Carleson measure for $B_p^0
(B)$, we need to prove
$$\sup\{{\mu(Q_\delta(\xi))\over {{\log^{1-p}{2\over \delta}
}}};\xi\in S, \delta>0\}<\infty.$$ For any $\xi\in S$, and
$0<\delta<1$, we consider the functions
$$f_{\xi,\delta}(z)=\log^{-{1\over
p}}\frac{2}{\delta}\log\frac{2}{1-\langle z,(1-\delta)\xi\rangle}.$$

Since $\log\frac{2}{1-\langle z,w\rangle}=\log 2+
\sum\limits_{k=1}^{\infty} k^{-1} {\langle z,w\rangle}^{k}$, and by
induction
 $$R^{\alpha,m}f(z)=\sum\limits_{k=0}^{\infty}\frac{(n+\alpha+1+k)\cdots (n+\alpha+m+k)}{(n+\alpha+1)\cdots (n+\alpha+m)}f_k (z),\ \ m\in N,$$
 we can get
$$\left|R^{\alpha,m}\log\frac{2}{1-\langle z,w\rangle}\right|\approx
\left|\sum\limits_{k=0}^{\infty} k^{m-1} {\langle
z,w\rangle}^{k}\right|\approx\left|{(1-\langle z,w\rangle
)}^{-m}\right |,$$ for all $w\in B$. The last formula is due to
${(1-\langle
z,w\rangle)}^{-m}=\sum\limits_{k=0}^{\infty}\frac{\Gamma (k+m)}{k!
\Gamma (m)} {\langle z,w\rangle}^k$.   Thus, using Proposition
1.4.10 of \cite{Ru} again, we know
\begin{eqnarray*}
\int_{Q_\delta (\xi)}{|f_{\xi,\delta}(z)|}^p d\mu
(z)&\leq&\int_{B}{|f_{\xi,\delta}(z)|}^p d\mu (z)\\
&\leq&C\int_B{|R^{\alpha,m} f_{\xi,\delta}(z)|}^p{(1-{|z|}^2)}^{mp-n-1}dv(z)\\
&\leq&C\int_B\frac{{(1-{|z|}^2)}^{mp-n-1}}{{|1-\langle z,(1-\delta)\xi\rangle|}^{mp}}\log^{-1}\frac{2}{\delta}dv(z)\\
&=&C\log^{-1}\frac{2}{\delta}\int_B\frac{{(1-{|z|}^2)}^{mp-n-1}}{{|
1-\langle
z,(1-\delta)\xi\rangle|}^{(n+1)+(mp-n-1)}}dv(z)\\
&\leq&C\log^{-1}\frac{2}{\delta}\log\frac{2}{\delta}\\
&\leq&C,
\end{eqnarray*}
where the condition $m+\sigma>{n\over p}$ is applied. By Lemma 2.6
of \cite{OF2}, we have $|f_{\xi,\delta}(z)|\approx\log^{1-{1\over
p}}\frac{2}{\delta}$ for $z\in Q_\delta (\xi)$. Consequently,
$$\sup\{{\mu(Q_\delta(\xi))\over {{\log^{1-p}{2\over \delta}
}}};\xi\in S, \delta>0\}<\infty.$$

On the other hand, suppose
$$\sup\{{\mu(Q_\delta(\xi))\over {{\log^{1-p-\varepsilon}{2\over
\delta} }}};\xi\in S, \delta>0\}<\infty,$$ using (2.3) provided
$M>mp$ large enough and by Lemma 2.3, we have
\begin{eqnarray*}
\int_B{|f(z)|}^p d\mu (z)&\leq&C\int_{B}{\left(\int_B | R^{\alpha,
m} f(w)|\frac{{(1-{|w|}^2)}^M}{{|1-\langle
z,w\rangle|}^{n+M+1-m}}dv(w)\right)}^p d\mu(z)\\
&\leq&C{\| f\|}_{B_p^0 (B)}^p,
\end{eqnarray*}
 this implies $\mu$ is a Carleson measure
for $B_p^0 (B)$.

{\bf Remark 2.1}\ \ Theorem 2.1 is an extension of Theorem 1 in
\cite{GP} to the higher dimensions. Since $p$-Carleson measure is
weaker than Carleson measures for $B_p^\sigma (B)$ for the case
$p=q$, it is natural that there exists a difference of an
arbitrarily small $\varepsilon>0$ between the necessity and the
sufficiency. However, we also note that such necessary and
sufficient conditions are unified, i.e. $\varepsilon=0$ for Hardy
spaces $H^p$ and the weighted Bergman spaces $A_\alpha^p$ in the
unit disk of ${\bf C}$ (see \cite{Car}, \cite{Du} and \cite{Lu}),
and the weighted Bergman spaces $A_\alpha^p$ in the unit ball of
${\bf C}^n$ (see \cite{Zh2}) in the case $p=q$.

\vspace{0.5cm}
\begin{center}
{\small\bf\S3\ \  Riemann-Stieltjes operators and
 multipliers for $B_p^\sigma (B)$ }
\end{center}
\vspace{0.3cm} \setcounter{section}{3} \setcounter{equation}{0}

\vspace{0.1cm}

In the following, we will apply these results to characterize the
Riemann-Stieltjes operators and
 multipliers for $B_p^\sigma (B)$ in the unit ball of ${\bf C}^n$.

\vspace{0.1cm}

{\bf Theorem 3.1}\ \ Suppose that $\varphi\in H(B)$, $1<p<\infty$,
$0\leq\sigma<\infty$, $m=1$, $1+\sigma>{n\over p}$. Then $U_\varphi:
B_p^\sigma (B)\rightarrow B_p^\sigma (B)$ is bounded if and only if
$\|\varphi\|_{H^\infty}<\infty$.

{\bf Proof}\ \ Noting that  $|R^{\alpha,1} f(z)|\approx |R f (z)|$,
we can work with the radial derivative $Rf(z)$. If
$\|\varphi\|_{H^\infty}<\infty$, then
\begin{eqnarray*}
\int_B{|R(U_\varphi f)(z)|}^p{(1-{|z|}^2)}^{p+\sigma p-n-1}
dv(z)&=&\int_B{|\varphi(z)|}^p{|R f(z)|}^p{(1-{|z|}^2)}^{p+\sigma
p-n-1}
dv(z)\\
&\leq&\|\varphi\|_{H^\infty}^p{\| f\|}_{B_p^\sigma (B)}^p.
\end{eqnarray*}
So, $U_\varphi: B_p^\sigma (B)\rightarrow B_p^\sigma (B)$ is
bounded.

Conversely, suppose $U_\varphi: B_p^\sigma (B)\rightarrow B_p^\sigma
(B)$ is bounded. For each $w\in B$ near to the boundary with
$|w|>{2\over 3}$. Choosing $f_w (z)=\frac{(1-{|w|}^2)}{{(1-\langle
z,w\rangle)}^{1+\sigma}}$. By Proposition 1.4.10 of \cite{Ru}, we
have

\begin{eqnarray*}
\int_B{|R f_w(z)|}^p{(1-{|z|}^2)}^{p+\sigma p-n-1}
dv(z)&=&\int_B\frac{{(1-{|w|}^2)}^p{|\langle
z,w\rangle|}^p{(1-{|z|}^2)}^{p+\sigma p-n-1}}{{|1-\langle
z,w\rangle|}^{(2+\sigma)p}}dv(z)\\
&\leq&{(1-{|w|}^2)}^p\int_B\frac{{(1-{|z|}^2)}^{p+\sigma
p-n-1}}{{|1-\langle
z,w\rangle|}^{(2+\sigma)p}}dv(z)\\
&\leq& C{(1-{|w|}^2)}^p{(1-{|w|}^2)}^{-p}\\
&\leq& C,\ \ \ \ \ \ \ \ \ \ \ \ \ \ \ \ \ \ \ \ \ \ \ \ \ \ \ \ \ \
\ \ \ \ \ \ \ \ \ \ \ \ \ \ \ \ \ \ (3.1)
\end{eqnarray*}
this implies $\sup\limits_{w\in B}{\|f_w\|}_{B_p^\sigma (B)}\leq C$.
It is well known that $$v(E(w,{1\over 2}))\approx
{(1-{|w|}^2)}^{n+1},\ \ \ 1-{|w|}^2\approx 1-{|z|}^2 \approx
|1-\langle z, w\rangle|\ \ \ \ {\rm for}\ \ \ z\in E(w,{1\over
2}).$$ Also note that for $z\in E(w,{1\over 2})$, we have
$$1-{| \varphi_w(z)|}^2=\frac{(1-{|w|}^2)(1-{|z|}^2)}{{|1-\langle
z,w\rangle|}^2}>\frac{3}{4}.$$ Thus $$1-|\langle z,w\rangle|\leq |
1-\langle z,w\rangle|<\frac{2}{\sqrt{3}}{(1-{|w|}^2)}^{1\over
2}(1-{|z|}^2)^{1\over 2}\leq \frac{2}{\sqrt{3}}{(1-{|w|}^2)}^{1\over
2}<\frac{2}{\sqrt{3}}\cdot\frac{\sqrt{5}}{3}=\frac{2\sqrt{15}}{9},$$
this implies $|\langle z,w\rangle|>1-\frac{2\sqrt{15}}{9}$. By the
$\mathscr{M}$-subharmonicity of ${|\varphi (w)|}^p$, we have
\begin{eqnarray*}
{|\varphi (w)|}^p&\leq& C\frac{1}{v(E(w,{1\over
2}))}\int_{E(w,{1\over 2})}{|\varphi (z)|}^p dv(z)\\
 &\leq& C\frac{1}{{(1-{|w|}^2)}^{n+1}}\int_{E(w,{1\over 2})}{|\varphi
(z)|}^pdv(z)\\
&\leq&C\int_{E(w,{1\over 2})}\frac{{(1-{|w|}^2)}^{p}}{{|1-\langle
z,w\rangle|}^{(2+\sigma)p}}{|\varphi
(z)|}^p{(1-{|z|}^2)}^{p+\sigma p-n-1}dv(z)\\
&\leq&C\int_{E(w,{1\over 2})}\frac{{|\langle
z,w\rangle|}^p{(1-{|w|}^2)}^{p}}{{|1-\langle
z,w\rangle|}^{(2+\sigma)p}}{|\varphi
(z)|}^p{(1-{|z|}^2)}^{p+\sigma p-n-1}dv(z)\\
&\leq&C\int_{B}{|\varphi (z)|}^p{| R
f_w(z)|}^p{(1-{|z|}^2)}^{p+\sigma p-n-1}dv(z)\\
&\leq& C{\|U_\varphi (f_w)\|}_{B_p^\sigma (B)}^p\leq
C{\|U_\varphi\|}^p {\|f_w\|}_{B_p^\sigma (B)}^p\leq C,
\end{eqnarray*}
and consequently, $|\varphi (w)|\leq C$ for $|w|>\frac{2}{3}$. By
maximum modulus principle, we have $|\varphi (w)|\leq C$ for $w\in
B$. Thus $\varphi\in H^\infty$.

\vspace{0.1cm}

{\bf Theorem 3.2}\ \ Suppose that $\varphi\in H(B)$, $1<p<\infty$,
$0\leq\sigma<\infty$, $m=1$,  $1+\sigma>{n\over p}$. Then the
following conditions are equivalent:

 (i)\ $V_\varphi: B_p^\sigma (B)\rightarrow B_p^\sigma (B)$ is bounded.

 (ii)\ The positive Borel measure $\mu_{\varphi}$ in $B$
 defined by $$d\mu_{\varphi}(z)={|R\varphi
 (z)|}^p{(1-{|z|}^2)}^{p+\sigma p-n-1}dv(z)$$ is a Carleson measure for
 $B_p^\sigma (B)$.

{\bf Proof}\ \ Note that $R(V_\varphi f)(z)=f(z)R\varphi (z)$.
Suppose  $V_\varphi: B_p^\sigma (B)\rightarrow B_p^\sigma (B)$ is
bounded. Then
\begin{eqnarray*}
\int_B{| f(z)|}^pd\mu_{\varphi}(z)&=&\int_B{| f(z)|}^p{|
R\varphi(z)|}^p{(1-{|z|}^2)}^{p+\sigma p-n-1}dv(z)\\
&=&\int_B{|
R(V_\varphi f)(z)|}^p{(1-{|z|}^2)}^{p+\sigma p-n-1}dv(z)\\
&\leq&{\| V_\varphi f\|}_{B_p^\sigma (B)}^p\leq{\| V_\varphi
\|}^p{\|f\|}_{B_p^\sigma (B)}^p \leq C{\|f\|}_{B_p^\sigma (B)}^p.
\end{eqnarray*}
So, $d\mu_{\varphi}$ is a Carleson measure for $B_p^\sigma (B)$.

Conversely, suppose $d\mu_{\varphi}(z)={|R\varphi
 (z)|}^p{(1-{|z|}^2)}^{p+\sigma p-n-1}dv(z)$ is a Carleson measure for
 $B_p^\sigma (B)$.
\begin{eqnarray*}
{\| V_\varphi f\|}_{B_p^\sigma (B)}^p&=&\int_B{|R(V_\varphi
f)(z)|}^p{(1-{|z|}^2)}^{p+\sigma p-n-1}
dv(z)\\
&=&\int_B{|f(z)|}^p{|R \varphi(z)|}^p{(1-{|z|}^2)}^{p+\sigma p-n-1}
dv(z)\\
&=&\int_B{|f(z)|}^pd\mu_{\varphi}(z)\\
&\leq& C{\| f\|}_{B_p^\sigma (B)}^p,
\end{eqnarray*}
this implies $V_\varphi: B_p^\sigma (B)\rightarrow B_p^\sigma (B)$
is bounded.

\vspace{0.1cm}

{\bf Corollary  3.1}\ \ Suppose that $\varphi\in H(B)$,
$1<p<\infty$, $0\leq\sigma<\infty$, $m=1$, $1+\sigma>{n\over p}$.
Then the following conditions are equivalent:

 (i)\ $M_\varphi: B_p^\sigma (B)\rightarrow B_p^\sigma (B)$ is bounded.

 (ii)\ $\varphi\in H^\infty$ and the positive Borel measure $\mu_{\varphi}$ in $B$
 defined by $$d\mu_{\varphi}(z)={|R\varphi
 (z)|}^p{(1-{|z|}^2)}^{p+\sigma p-n-1}dv(z)$$ is a Carleson measure for
 $B_p^\sigma (B)$.

{\bf Proof}\ \ The implication $"(ii)\Rightarrow (i)"$ follows from
Theorem 3.1 and Theorem 3.2, using the fact that $R(M_\varphi
f)(z)=R(V_\varphi f)(z)+R(U_\varphi f)(z)$.

$"(i)\Rightarrow (ii)".$

Suppose $M_\varphi: B_p^\sigma (B)\rightarrow B_p^\sigma (B)$ is
bounded. At first, similar to the proof of Theorem 3.1, we will
prove that $\varphi\in H^\infty$.  For each $w\in B$ near to the
boundary with $|w|>{2\over 3}$, set
$$\widetilde{f_w}(z)=\frac{(1-{|w|}^2)}{{(1-\langle
z,w\rangle)}^{1+\sigma}}-{(1-{|w|}^2)}^{-\sigma},\ \ \ \ \ \ z\in
B.$$ Since $R\widetilde{f_w}(z)=R f_w(z)$, (3.1) implies that
$\sup\limits_{w\in B}{\|\widetilde{f_w}\|}_{B_p^\sigma (B)}\leq C$.
Noting that $\widetilde{f_w}(w)=0$, by the
$\mathscr{M}$-subharmonicity of ${|R(M_\varphi
\widetilde{f_w})(w)|}^p$, we have
\begin{eqnarray*}
{|\varphi (w)|}^p&\leq&
C{|w|}^{2p}{(1-{|w|}^2)}^{-(1+\sigma)p}{|\varphi
(w)|}^p{(1-{|w|}^2)}^{p+\sigma p}\\
&=&C{|R\widetilde{f_w}(w)|}^p{|\varphi
(w)|}^p{(1-{|w|}^2)}^{p+\sigma p}\\
&=&C{|R\widetilde{f_w}(w)\varphi(w)+\widetilde{f_w}(w)R\varphi(w)|}^p{(1-{|w|}^2)}^{p+\sigma p}\\
&=&C{|R(M_\varphi
\widetilde{f_w})(w)|}^p{(1-{|w|}^2)}^{p+\sigma p}\\
&\leq&C\frac{{(1-{|w|}^2)}^{p+\sigma p}}{v(E(w,{1\over
2}))}\int_{E(w,{1\over 2})}{|R(M_\varphi
\widetilde{f_w})(z)|}^pdv(z)\\
&\leq&C\int_{B}{|R(M_\varphi
\widetilde{f_w})(z)|}^p{(1-{|z|}^2)}^{p+\sigma p-n-1}dv(z)\\
&\leq&C{\| M_\varphi \widetilde{f_w}\|}_{B_p^\sigma (B)}^p\leq C{\|
M_\varphi \|}^p{\|\widetilde{f_w}\|}_{B_p^\sigma (B)}^p \leq C,
\end{eqnarray*}
and consequently, $|\varphi (w)|\leq C$ for $|w|>\frac{2}{3}$. By
maximum modulus principle, we have $|\varphi (w)|\leq C$ for $w\in
B$. Thus $\varphi\in H^\infty$.

{\bf Remark 3.1}\ As to the Riemann-Stieltjes operators and
multipliers on Besov-Sobolev  spaces, in the case of one complex
variable, there are a lot of results, see \cite{GP}, \cite{KS},
\cite{St},  \cite{Wu}. In the case of several complex variables, we
can find the research has been developing, see \cite{ARS1},
\cite{ARS2}, \cite{OF1}, \cite{PO2}. Such question on other spaces
was studied in \cite{Hu}, \cite{LS}, \cite{OF2}, \cite{OF4},
\cite{Xi1}, \cite{Xi2}.

\end{document}